\documentclass[12pt]{article}

\makeatletter

\def\eqnarray{\stepcounter{equation}\let\@currentlabel=\theequation
\global\@eqnswtrue
\global\@eqcnt\z@\tabskip\@centering\let\\=\@eqncr
$$\halign to \displaywidth\bgroup\@eqnsel\hskip\@centering
  $\displaystyle\tabskip\z@{##}$&\global\@eqcnt\@ne 
  \hfil$\displaystyle{\hbox{}##\hbox{}}$\hfil
  &\global\@eqcnt\tw@ $\displaystyle\tabskip\z@
  {##}$\hfil\tabskip\@centering&\llap{##}\tabskip\z@\cr}

\makeatother

\newbox\squ  
\setbox\squ=\hbox{\vrule width.3pt
            \vbox{\hrule height.3pt width.4em\kern1ex\hrule height.3pt}%
            \vrule width.3pt}
\def\endproof{%
  \ifmmode\eqno\copy\squ\medskip\else{\unskip\nobreak\hfil%
    \penalty50\hskip2em\hbox{}\nobreak\hfil\copy\squ
    \parfillskip=0pt \finalhyphendemerits=0\penalty-100\medskip}
  \fi} 
  
\usepackage{amssymb}
\usepackage{latexsym}

\begin{document}


\newcommand{\scl}{\scriptstyle}
\newcommand{\non}{\nonumber}
\newcommand{\bib}{thebibliography}
\newcommand{\tr}{ {\rm tr}\ts}
\newcommand{\wt}{\widetilde}
\newcommand{\wh}{\widehat}
\newcommand{\ot}{\otimes}
\newcommand{\g}{\mathfrak g}
\newcommand{\ts}{\,}
\newcommand{\Sym}{\mathfrak S}
\newcommand{\mod}{ {\rm mod}\ }
\newcommand{\middle}{ {\rm middle} }
\newcommand{\F}{ {\rm F}}
\newcommand{\A}{ {\rm A}}
\renewcommand{\P}{ {\rm P}}
\newcommand{\E}{\cal E}
\newcommand{\U}{ {\rm U}}
\newcommand{\Z}{ {\rm Z}}
\newcommand{\Norm}{ {\rm Norm}\ts}
\newcommand{\R}{ {\rm R}}
\newcommand{\I}{ {\rm I}}
\newcommand{\J}{ {\rm J}}
\newcommand{\Y}{ {\rm Y}}
\newcommand{\Srm}{ {\rm S}}
\newcommand{\End}{ {\rm End}\ts}
\newcommand{\C}{\mathbb{C}}
\newcommand{\ZZ}{\mathbb{Z}}
\newcommand{\dega}{ {\rm deg}\ts}
\newcommand{\gr}{ {\rm gr}\ts}
\newcommand{\sgn}{ {\rm sgn}}
\newcommand{\qdet}{ {\rm qdet}\ts}
\newcommand{\tra}{ {\rm tr}\ts}
\newcommand{\T}{{\rm T}}
\newcommand{\agot}{\mathfrak a}
\newcommand{\h}{\mathfrak h}
\newcommand{\n}{\mathfrak n}
\newcommand{\sll}{\mathfrak{sl}}
\newcommand{\gl}{\mathfrak{gl}}
\newcommand{\Proof}{\noindent{\bf Proof.}\ \ }   
\renewcommand{\theequation}{\arabic{section}.\arabic{equation}}  

\newcommand{\ble}{\begin{lem}}
\newcommand{\ele}{\end{lem}}
\newcommand{\bth}{\begin{thm}}
\renewcommand{\eth}{\end{thm}}
\newcommand{\bpr}{\begin{prop}}
\newcommand{\epr}{\end{prop}}
\newcommand{\bco}{\begin{cor}}
\newcommand{\eco}{\end{cor}}
\newcommand{\bde}{\begin{defin}}
\newcommand{\ede}{\end{defin}}

\newcommand{\beq}{\begin{equation}}
\newcommand{\eeq}{\end{equation}}
\newcommand{\bea}{\begin{eqnarray}}
\newcommand{\eea}{\end{eqnarray}}

\newtheorem{thm}{Theorem}[section]
\newtheorem{prop}[thm]{Proposition}
\newtheorem{cor}[thm]{Corollary}
\newtheorem{defin}[thm]{Definition}
\newtheorem{lem}[thm]{Lemma}

\title{\LARGE\bf Yangians and transvector algebras}
\author{{\sc A. I. Molev}\\[15pt]
School of Mathematics and Statistics\\
University of Sydney,
NSW 2006, Australia\\
alexm@maths.usyd.edu.au\\[30pt]
Research Report 98--30}

\date{November 1998}
\maketitle

\vspace{7 mm}
\begin{abstract}
Olshanski'\,s centralizer construction
provides a realization of the Yangian $\Y(m)$
for the Lie algebra $\gl(m)$ as a subalgebra in the projective limit algebra 
$\A_m={\rm lim\ts proj\ts}\A_m(n)$ 
as $n\to\infty$, where $\A_m(n)$ is the centralizer
of $\gl(n-m)$
in the enveloping algebra $\U(\gl(n))$.
We give a modified version
of this construction based on a quantum
analog of Sylvester's theorem. 
We then use it to get an algebra
homomorphism from the Yangian $\Y(n)$
to the transvector algebra
associated with the pair $\gl(m)\subset\gl(m+n)$.
The results are applied to identify
the elementary representations of the Yangian
by constructing their highest vectors
explicitly in terms of elements of
the transvector algebra.
\end{abstract}


\pagebreak

\setcounter{section}{-1}

\section{Introduction}   \label{i}
\setcounter{equation}{0}

Let $\A_m(n)$ denote the centralizer
of $\gl(n-m)$ in the universal enveloping algebra $\U(\gl(n))$.
In particular, $\A_0(n)$ is the center of $\U(\gl(n))$.
It was shown by Olshanski~\cite{o:ea,o:ri} 
that for any fixed $m$
there exists a chain of natural homomorphisms
\bea\label{chainA}
\A_m(m)\leftarrow\A_m(m+1)
\leftarrow\cdots
\leftarrow\A_m(n)
\leftarrow\cdots
\eea
and one can define the corresponding
projective limit algebra
\bea   \label{lim}
\A_m={\rm lim\ts proj\ts} \A_m(n),
\qquad n\to\infty.
\eea
The algebra $\A_0$ is
isomorphic to an algebra of polynomials in countably many
variables while for $m>0$ one has the tensor 
product decomposition~\cite{o:ea,o:ri}:
\bea  \label{tensor}
\A_m=\A_0\otimes \Y(m),
\eea
where $\Y(m)=\Y(\gl(m))$ is the {\it Yangian\/} for the Lie algebra $\gl(m)$;
see the definition in Section~\ref{qst} below.
The algebra $\Y(m)$
first appeared in the works of
Faddeev's school on the Yang--Baxter equation; 
see e.g.~\cite{krs:yb,ks:qs,tf:qi}. 
The Yangian $\Y(m)$ is a Hopf algebra
which can be regarded as a deformation of the enveloping algebra
$\U(\gl(m)[x])$, where $\gl(m)[x]$
is the Lie algebra of $\gl(m)$-valued polynomials~\cite{d:ha}.

The key part of the proof of the
decomposition~(\ref{tensor}) is a construction
of algebra homomorphisms
\bea   \label{hom}
\Y(m)\to\A_m(n),\qquad n=m,m+1,\dots
\eea
compatible with the chain (\ref{chainA}). Then one shows that
this defines an embedding
\bea   \label{embed}
\Y(m)\hookrightarrow \A_m
\eea
so that the Yangian $\Y(m)$ can be identified 
with a subalgebra in $\A_m$; see~\cite{o:ea,o:ri}. 

Note that
a projective limit algebra of type (\ref{lim}) can be also constructed
for the series of the orthogonal and
symplectic Lie algebras. This limit algebra turns out to be related with
the corresponding twisted Yangians; see \cite{mo:cc,o:ty} for more
details.

In this paper we construct a new family of
homomorphisms~(\ref{hom}) which define 
a new embedding
(\ref{embed}).
Then we give a modified tensor product decomposition of type~(\ref{tensor}).
Our argument is based on a quantum analog of Sylvester's theorem
for the Yangians (Theorem~\ref{thm:qs}); cf. \cite{gr:dm, kl:mi}.

We also use Theorem~\ref{thm:qs}
to identify the so-called elementary representations
of $\Y(n)$. They naturally arise
from the centralizer construction as follows. 
Consider
a finite-dimensional irreducible representation $L(\lambda)$ of the
Lie algebra $\gl(m+n)$ with the highest weight $\lambda$
and denote by $L(\lambda)^+_{\mu}$ the subspace in $L(\lambda)$ of
$\gl(m)$-highest vectors of weight $\mu$. 
It is well-known (see e.g.~\cite[Section~9.1]{d:ae}) that 
$L(\lambda)^+_{\mu}$ is an irreducible representation of the
centralizer $\A=\U(\gl(m+n))^{\gl(m)}$. We can make $L(\lambda)^+_{\mu}$
into a $\Y(n)$-module via the homomorphism
$\Y(n)\to \A$; cf.~(\ref{hom}). This module can be shown to be irreducible;
see Section~\ref{hvY}. 

On the other hand, the subspace $L(\lambda)^+$
of $\gl(m)$-highest vectors in $L(\lambda)$
is preserved by the so-called {\it raising\/} and 
{\it lowering operators\/}. The latter generate 
the {\it transvector algebra\/} $\Z$ 
associated with the pair $\gl(m)\subset \gl(m+n)$
(which is sometimes called the $S$-{\it algebra\/} 
or {\it Mickelsson algebra\/});
see~\cite{z:sa}--\cite{z:it}.
By definition of $\Z$ (see Section~\ref{tc})
there is a natural algebra homomorphism $\A\to\Z$ so that
the action of generators of $\Y(n)$ in $L(\lambda)^+_{\mu}$
can be explicitly expressed in terms 
of the lowering and raising operators (see Theorem~\ref{thm:images}).

We explicitly construct the highest vector
of the $\Y(n)$-module $L(\lambda)^+_{\mu}$
in terms of the lowering operators
and calculate its highest weight.
We also identify $L(\lambda)^+_{\mu}$ 
as a $\Y(\sll(n))$-module by calculating its Drinfeld
polynomials. They turn out to be the same as those
found by Nazarov and Tarasov~\cite{nt:ry} for
a different representation of $\Y(\sll(n))$ in $L(\lambda)^+_{\mu}$
defined via the Olshanski homomorphism (\ref{hom}),
which shows that these two $\Y(\sll(n))$-modules are isomorphic.

Following \cite{nt:ry} we call the $\Y(n)$-module $L(\lambda)^+_{\mu}$
{\it elementary\/}. These representations play an important role in the
classification of the $\Y(n)$-modules with a semisimple
action of the Gelfand--Tsetlin subalgebra; see~\cite{c:ni,nt:ry}.
In particular,
it was proved in~\cite[Theorem~4.1]{nt:ry} that, up to
an automorphism of $\Y(n)$, any such
module is isomorphic to a tensor product of
elementary representations.

\section{Quantum Sylvester's theorem}\label{qst}
\setcounter{equation}{0}

A detailed description of the algebraic structure
of the Yangian for the Lie algebra $\gl(n)$ is given
in the expository paper~\cite{mno:yc}. In this section
we reproduce some of those results and use them to prove
a quantum analog of Sylvester's theorem.
 
The {\it Yangian\/} $\Y(n)=\Y(\gl(n))$ is the
complex associative algebra with the
generators $t_{ij}^{(1)},t_{ij}^{(2)},\dots$ where $1\leq i,j\leq n$,
and the defining relations
\bea\label{defrel}
[t_{ij}(u),t_{kl}(v)]=\frac{1}{u-v}(t_{kj}(u)t_{il}(v)-t_{kj}(v)t_{il}(u)),
\eea
where
\bea\label{series}
t_{ij}(u) = \delta_{ij} + t^{(1)}_{ij} u^{-1} + t^{(2)}_{ij}u^{-2} +
\cdots \in \Y(n)[[u^{-1}]]
\non
\eea
and $u$ is a formal variable.
Introduce the matrix
\bea\label{tmatrix}
T(u):=\sum_{i,j=1}^n t_{ij}(u)\ot E_{ij}
\in \Y(n)[[u^{-1}]]\ot \End\C^n,
\non
\eea
where the $E_{ij}$ are the standard matrix units.
Then the relations~(\ref{defrel}) are equivalent to
the single relation
\bea\label{ternary}
R(u-v)T_1(u)T_2 (v) = T_2(v)T_1(u)R(u-v).
\eea
Here $T_1(u)$ and $T_2(u)$ are regarded as elements of
$\Y(n)[[u^{-1}]]\ot \End\C^n\ot \End\C^n$, the subindex of $T(u)$
indicates to which copy of $\End\C^n$ this matrix corresponds, and
\bea\label{R(u)}
R(u)=1-Pu^{-1},\qquad P=\sum_{i,j=1}^n E_{ij}\ot E_{ji}\in
(\End\C^n)^{\otimes 2}.
\non
\eea

The {\it quantum determinant\/} $\qdet T(u)$
of the matrix $T(u)$ is a formal series in $u^{-1}$
with coefficients from $\Y(n)$ defined by
\bea\label{qdety}
\qdet T(u)=\sum_{p\in \Sym_n} \sgn(p)\ts t_{p(1)1}(u)\cdots
t_{p(n)n}(u-n+1).
\eea

The coefficients of the quantum
determinant $\qdet T(u)$ are algebraically independent generators of
the center of the algebra $\Y(n)$.

We shall need a generalization of the relation~(\ref{ternary})
for elements of multiple tensor products of the form
$\Y(n)[[u^{-1}]]\ot \End\C^n\ot\cdots\ot \End\C^n$.
For an operator $X\in \End\C^n$
and a number $s=1,2,\ldots$
we set 
\bea\label{X_i}
X_i:=1^{\ot (i-1)}\ot X\ot 1^{\ot (s-i)}
\in(\End\C^n)^{\ot s}, \quad 1\leq i\leq s.
\non
\eea
Similarly,
if $X\in(\End\C^n)^{\otimes 2}$ then for any $i,j$\  
such that 
$1\leq i,j\leq s$\ and $i\neq j$, we denote by $X_{ij}$ the operator in 
$(\C^n)^{\ot s}$
which acts as $X$ in the product of $i$th and $j$th copies and as
1 in all other copies. 
Let $u_1,\ldots,u_s$ be formal variables. Set
\bea\label{train}
R(u_1,\dots ,u_s) :=(R_{s-1,s})(R_{s-2,s}R_{s-2,s-1}) \cdots (R_{1s}
\cdots R_{12}),
\eea
where  we abbreviate $R_{ij}:=R_{ij}(u_i-u_j)$.
Note that the following (Yang--Baxter) relation
is satisfied by the $R_{ij}$
\bea\label{ybe}
R_{ij}R_{ir}R_{jr}=R_{jr}R_{ir}R_{ij},
\eea
and the factors $R_{ij}$ and $R_{rs}$ with distinct indices
are permutable in~(\ref{train}). This allows one to deduce
the following identity from (\ref{ternary}):
\bea\label{fundam}
R(u_1,\dots, u_s)\ts T_1(u_1) \cdots T_s(u_s) = 
T_s(u_s) \cdots T_1(u_1)\ts R(u_1,\dots, u_s).
\eea

The symmetric group $\Sym_s$ naturally
acts in the tensor space $(\C^n)^{\ot s}$
by permutations of the tensor factors. If the variables in~(\ref{train})
satisfy the condition $u_i-u_{i+1}=1$ for $i=1,\dots,s-1$ then
$R(u_1,\dots ,u_s)$ becomes the antisymmetrization operator in
$(\C^n)^{\ot s}$:
\bea\label{anti}
R(u_1,\dots ,u_s)=A_s:=\sum_{q\in \Sym_s}\sgn(q)\ts Q,
\eea
where $Q$ is the operator in $(\C^n)^{\ot s}$ corresponding
to a permutation $q\in \Sym_s$.
By (\ref{fundam}) we then have
\beq\label{funA}
A_s\ts T_1(u) \cdots T_s(u-s+1) = 
T_s(u-s+1) \cdots T_1(u)\ts A_s.
\eeq
The operator (\ref{funA}) can be written in the form
\bea\label{matelem}
\sum {t\ts}^{a_1\cdots\ts a_s}_{b_1\cdots\ts b_s}(u)\ot E_{a_1b_1}\ot \cdots
\ot E_{a_sb_s},
\non
\eea
summed over the indices $a_i,b_i\in\{1,\dots,n\}$, where 
${t\ts}^{a_1\cdots\ts a_s}_{b_1\cdots\ts b_s}(u)\in \Y(n)[[u^{-1}]]$.
In particular,
${t\ts}^{a}_{b}(u)=t_{ab}(u)$. Note that by~(\ref{funA})
the series ${t\ts}^{a_1\cdots\ts a_s}_{b_1\cdots\ts b_s}(u)$ is antisymmetric
with respect to permutations of the upper indices and of
the lower indices. It can be
given by the following explicit formulas
\bea\label{qminor}
{t\ts}^{a_1\cdots\ts a_s}_{b_1\cdots\ts b_s}(u)&=&
\sum_{\sigma\in \Sym_s} \sgn(\sigma)\ts t_{a_{\sigma(1)}b_1}(u)\cdots
t_{a_{\sigma(s)}b_s}(u-s+1)\\
\label{qminor2}
\mbox{}&=&
\sum_{\sigma\in \Sym_s} \sgn(\sigma)\ts t_{a_1b_{\sigma(1)}}(u-s+1)\cdots
t_{a_sb_{\sigma(s)}}(u).
\eea
Note that ${t\ts}^{1\cdots\ts n}_{1\cdots\ts n}(u)=\qdet T(u)$;
see (\ref{qdety}).

\bpr\label{prop:qmr}
We have the relations
\bea\label{qmrel}
\lefteqn{
[{t\ts}^{a_1\cdots\ts a_k}_{b_1\cdots\ts b_k}(u),
{t\ts}^{c_1\cdots\ts c_l}_{d_1\cdots\ts d_l}(v)]=
\sum_{p=1}^{\min\{k,l\}}\frac{(-1)^{p-1}\ts p!}
{(u-v-k+1)\cdots (u-v-k+p)}
}
\non
\\
&&\sum_{\scl i_1<\cdots<i_p\atop
\scl j_1<\cdots<j_p}\left(
{t\ts}^{a_1\cdots\ts c_{j_1}\cdots\ts 
c_{j_p}\cdots\ts a_k}_{b_1\ \cdots\ \  b_k}(u)
{t\ts}^{c_1\cdots\ts a_{i_1}\cdots\ts 
a_{i_p}\cdots\ts c_l}_{d_1\ \cdots\ \  d_l}(v)
-{t\ts}^{c_1\ \cdots\ \  c_l}_{d_1\cdots\ts 
b_{i_1}\cdots\ts b_{i_p}\cdots\ts d_l}(v)
{t\ts}^{a_1\ \cdots\ \  a_k}_{b_1\cdots\ts 
d_{j_1}\cdots\ts d_{j_p}\cdots\ts b_k}(u)
\right).
\non
\eea
Here the $p$-tuples
of upper indices $(a_{i_1},\dots, a_{i_p})$ and $(c_{j_1},\dots, c_{j_p})$
are respectively
interchanged in the first summand on the right hand
side while the $p$-tuples of lower indices
$(b_{i_1},\dots, b_{i_p})$ and
$(d_{j_1},\dots, d_{j_p})$
are interchanged in the second
summand.
\epr

\Proof Let us take $s=k+l$ in~(\ref{fundam}) and specialize
the variables $u_i$ as follows:
\beq\label{special}
u_i=u-i+1,\quad i=1,\dots,k \quad{\rm and}\quad
u_{k+j}=v-j+1,\quad j=1,\dots,l.
\eeq
Using (\ref{ybe}) we can bring the product
(\ref{train}) to the form
\beq\label{R=}
R(u_1,\dots, u_{k+l})=\left(\prod_{j=1,\dots,\ts l}^{\to} 
R_{k,k+j}\cdots R_{1,k+j}\right)R(u_1,\dots, u_{k})\ts 
R(u_{k+1},\dots, u_{k+l}).
\eeq
However, by (\ref{anti}) and (\ref{special})
\beq\label{RA}
R(u_1,\dots, u_{k})=A_k,\qquad 
R(u_{k+1},\dots, u_{k+l})=A'_l,
\eeq
where $A'_l$ is the antisymmetrizer 
corresponding to the set of
indices $\{k+1,\dots,k+l\}$. Let us show that the expression~(\ref{R=})
can then be written as
\beq\label{Rti}
\wt R(u,v)\ts A_kA'_l
\eeq
with
\beq\label{Rti=}
\wt R(u,v)=
\sum_{p=0}^{\min\{k,l\}}\frac{(-1)^{p}\ts p!}
{(u-v-k+1)\cdots (u-v-k+p)}
\sum_{\scl i_1<\cdots<i_p\atop
\scl j_1<\cdots<j_p}P_{i_1,k+j_1}\cdots P_{i_p,k+j_p}.
\eeq
Indeed, if indices $i_1,\dots,i_s\in\{1,\dots,k\}$ are distinct then
\bea\label{PPA}
P_{i_1,k+j}\cdots P_{i_s,k+j}\ts A_k=(-1)^{s-1} P_{i_s,k+j}\ts A_k.
\non
\eea
This implies that for each $j$
\bea\label{RRA_k}
R_{k,k+j}\cdots R_{1,k+j}\ts A_k=
\left(1-\frac{P_{1,k+j}+\cdots+P_{k,k+j}}{u-v-k+j}\right) A_k.
\eea
Similarly, if indices $j_1,\dots,j_s\in\{1,\dots,l\}$ are distinct then
\bea\label{PPA'}
P_{i,k+j_1}\cdots P_{i,k+j_s}\ts A'_l=(-1)^{s-1} P_{i,k+j_s}\ts A'_l.
\eea
Moreover, if $i_1\ne i_2$ and $j_1\ne j_2$ then
\bea\label{iijj}
P_{i_1,k+j_1}P_{i_2,k+j_2}\ts A_kA'_l=P_{i_1,k+j_2}P_{i_2,k+j_1}\ts A_kA'_l.
\non
\eea
Using this relation together with~(\ref{PPA'}) we bring~(\ref{R=})
to the form~(\ref{Rti}) by a straightforward calculation.
Now~(\ref{fundam}) takes the form
\bea
\wt R(u,v)A_kT_1(u)&\cdots& T_k(u-k+1)\ts 
A'_l\ts T_{k+1}(v)\cdots T_{k+l}(v-l+1)
\non
\\
=A'_l\ts T_{k+1}(v)&\cdots& T_{k+l}(v-l+1)\ts A_kT_1(u)\cdots T_k(u-k+1)
\wt R(u,v),\label{RATAT}
\eea
where we have used the identity
(\ref{funA}) and the fact that $\wt R(u,v)\ts A_kA'_l=A_kA'_l\ts \wt R(u,v)$
which is easy to verify.
To complete the proof we take the coefficients at
$E_{a_1b_1}\ot\cdots\ot E_{a_kb_k}\ot E_{c_1d_1}\ot\cdots\ot E_{c_ld_l}$
on both sides of (\ref{RATAT})
(or, equivalently, apply both sides to the
basis vector $e_{b_1}\ot\cdots\ot e_{b_k}\ot e_{d_1} \ot\cdots\ot e_{d_l}$
of $(\C^n)^{\ot (k+l)}$ and compare the coefficients at the vector
$e_{a_1}\ot\cdots\ot e_{a_k}\ot e_{c_1} \ot\cdots\ot e_{c_l}$, where
$e_1,\dots,e_n$ is the canonical basis of $\C^n$).
\endproof


Introduce the {\it quantum comatrix\/}
$\wh T(u)=(\wh t_{ij}(u))$ for the matrix $T(u)$ by
\bea\label{qcom}
\wh t_{ij}(u)=(-1)^{i+j}
{t\ts}^{1\ts\cdots\ts
\wh j\ts\cdots\ts n}_{1\ts\cdots\ts\wh i\ts\cdots\ts n}(u),
\non
\eea
where the hats on the right hand side
indicate the indices to be omitted.
Equivalently, $\wh T(u)$ can be defined by
\bea\label{Th}
A_nT_1(u)\cdots T_{n-1}(u-n+2)=A_n\wh T_n(u).
\non
\eea
Taking $s=n$ in (\ref{funA}) we derive
the identity
\bea\label{t:hatt}
\wh T(u)T(u-n+1)=\qdet T(u).
\eea

The Poincar\'e--Birkhoff--Witt theorem for the algebra $\Y(n)$
(see e.g. \cite[Corollary~1.23]{mno:yc}) implies that
for $m\leq n$ the Yangians $\Y(m)$ and $\Y(n-m)$ can be identified
with the subalgebras in $\Y(n)$ generated by the coefficients
of the series $t_{ij}(u)$ with $i,j\leq m$ and 
$m<i,j\leq n$, respectively.
For any indices $i,j\leq m$ introduce the following series
with coefficients in $\Y(n)$
\bea\label{ttilde}
{\wt t}_{ij}(u)={t\ts}^{i,m+1\cdots\ts n}_{j,m+1\cdots\ts n}(u)
\eea
and combine them into the matrix $\wt T(u)=(\wt t_{ij}(u))$.
For subsets ${\cal P}$ and ${\cal Q}$ in $\{1,\dots,n\}$
and an $n\times n$-matrix $X$
we shall denote by $X^{}_{{\cal P}{\cal Q}}$ the submatrix
of $X$ whose rows and columns are enumerated by 
${\cal P}$ and ${\cal Q}$ respectively.
Set ${\cal A}=\{1,\dots,m\}$ and ${\cal B}=\{m+1,\dots,n\}$.
We shall need the following generalization of (\ref{t:hatt}).

\bpr\label{prop:tAA}
We have the identity
\bea\label{T:hatAA}
\wh T(u)^{}_{{\cal A}{\cal A}}\ts\wt T(u-m+1)=\qdet T(u)\ts\qdet
T(u-m+1)^{}_{{\cal B}{\cal B}}.
\eea
\epr

\Proof For $m=1$ the relation (\ref{T:hatAA}) is trivial.
We shall assume 
that $m\geq 2$.
Consider
the identity (\ref{fundam}) with $s=2n-m$
and specialize the variables $u_i$ as follows
\[
u_i=u-i+1,\quad i=1,\dots,n-1 \quad{\rm and}\quad
u_{n+j-1}=v-j+1,\quad j=1,\dots,n-m+1.
\]
We shall show that (\ref{fundam}) then becomes a rational function
in $v$ with a simple pole at $v=u-m+1$ and
calculate the corresponding residue in two different ways. First,
write (\ref{fundam}) in the form (\ref{RATAT})
with $k=n-1$ and $l=n-m+1$. Then 
fix indices $i,j\in\{1,\dots,m\}$ and
apply the left hand side of (\ref{RATAT})
to the basis vector
\bea\label{vecbeg}
e_1\ot\cdots\ot e_{i-1}\ot e_{i+1} \ot 
\cdots\ot e_n\ot e_j\ot e_{m+1} \ot 
\cdots\ot e_n
\eea
in $(\C^n)^{\ot(2n-m)}$ and take the coefficient at the vector
\bea\label{vecend}
e_1\ot\cdots\ot e_{i-1}\ot e_{i+1} \ot 
\cdots\ot e_n\ot e_i\ot e_{m+1} \ot 
\cdots\ot e_n.
\eea
A straightforward calculation shows that
this matrix element equals
\bea\label{matel}
\frac{u-v-n+1}{u-v-m+1}\left(\wh t_{ii}(u)\wt t_{ij}(v)
+\frac{1}{u-v-m+2}\sum_{a=1,\ts a\ne i}^{m}\wh t_{ia}(u)\wt t_{aj}(v)\right).
\non
\eea
Multiplying this expression by $u-v-m+1$ and setting $v=u-m+1$
we obtain
\bea\label{lefths}
(m-n)\left(\wh T(u)^{}_{{\cal A}{\cal A}}\ts\wt T(u-m+1)\right)_{ij}.
\non
\eea

On the other hand, transform $R(u_1,\dots,u_{2n-m})$ applying
(\ref{R=}) and (\ref{RA}) with $k=n-1$ and $l=n-m+1$. 
Using (\ref{RRA_k}) we can write $R(u_1,\dots,u_{2n-m})$ as
\bea\label{RAAR}
\left(\prod_{j=1,\dots,\ts l-1}^{\to} 
R_{k,k+j}\cdots R_{1,k+j}\right)
\left(1-\frac{P_{1,k+l}+\cdots+P_{k,k+l}}{u-v-k+l}\right)A_kA'_l.
\eea
Note that $(l-1)!\ts A'_l=A'_{l-1}A'_l$. Therefore, 
repeating the argument of the proof of Proposition~\ref{prop:qmr}
we can transform (\ref{RAAR}) to get the expression
\beq\label{Rtil}
\frac{1}{(l-1)!}\ts\wt R(u,v)A_kA'_{l-1}
\left(1-\frac{P_{1,k+l}+\cdots+P_{k,k+l}}{u-v-k+l}\right)
A'_l,
\eeq
where $\wt R(u,v)$ is given by (\ref{Rti=}) with $l$ replaced by $l-1$.
Finally, put $k=n-1$, $l=n-m+1$ into (\ref{Rtil}),
multiply it by $u-v-m+1$ and then put $v=u-m+1$. This gives
\beq\label{PPAA}
(m-n)\sum
P_{i_1,n}\cdots P_{i_{n-m},2n-m-1}
(1-P_{1,2n-m}-\cdots-P_{n-1,2n-m})A_{n-1}A'_{n-m+1},
\eeq
summed over the sets of indices $1\leq i_1<\cdots<i_{n-m}\leq n-1$.
Note that 
\bea\label{anti+}
\wt A_n:=(1-P_{1,2n-m}-\cdots-P_{n-1,2n-m})A_{n-1}
\non
\eea
is the antisymmetrizer corresponding
to the subset of indices $\{1,\dots,n-1,2n-m\}$.
The application of the operator (\ref{PPAA})
to the vector (\ref{vecbeg}) gives the zero vector unless $i=j$
since it is annihilated by $\wt A_nA'_{n-m+1}$. Suppose now that
$i=j$. By (\ref{funA}) we can write
\bea\label{ATT=}
A'_{n-m+1}T_{n}(u-m+1)\cdots &\ts\displaystyle{T_{2n-m}(u-n+1)}\ts &=
\non
\\
&\displaystyle{T_{2n-m}(u-n+1)}&\cdots \ts T_{n}(u-m+1) A'_{n-m+1}.
\non
\eea
The residue of the left hand side of (\ref{fundam})
at $v=u-m+1$ then
takes the form
\bea
(m-n)\sum P_{i_1,n}\cdots 
P_{i_{n-m},2n-m-1}
\wt A_n T_1(u)&\cdots & T_{n-1}(u-n+2)T_{2n-m}(u-n+1)
\non
\\
T_{2n-m-1}(u-n+2)&\cdots &T_{n}(u-m+1)A'_{n-m+1}.\label{fundamend}
\non
\eea
The diagonal matrix element of this operator corresponding to
the vector (\ref{vecend})
equals
\bea\label{rightend}
(m-n)\ts \qdet T(u)\ts\qdet
T(u-m+1)^{}_{{\cal B}{\cal B}}
\non
\eea
which completes the proof. \endproof


We are now in a position to prove a quantum analog
of Sylvester's theorem for the Yangians (Theorem~\ref{thm:qs}). 
An analog of Sylvester's theorem for
the quantized algebra of functions on $GL(n)$ 
was given
by Krob and Leclerc~\cite{kl:mi} with the use 
of the quasi-determinant version
of this theorem due to Gelfand and Retakh~\cite{gr:dm}.
As was remarked in \cite{kl:mi}, 
the same approach can be used for the case of the Yangians.
We give a different proof based on 
Propositions~\ref{prop:qmr} and \ref{prop:tAA}.
We use notation (\ref{ttilde}).

\bth\label{thm:qs}
The mapping
\bea\label{homtt}
t_{ij}(u)\mapsto {\wt t}_{ij}(u),\qquad 1\leq i,j\leq m
\eea
defines an algebra homomorphism $\Y(m)\to\Y(n)$. Moreover,
one has the identity
\bea\label{qdettt}
\qdet {\wt T}(u)=\qdet T(u) \ts
\qdet T(u-1)^{}_{{\cal B}{\cal B}}\cdots\qdet
T(u-m+1)^{}_{{\cal B}{\cal B}}. 
\eea
\eth

\Proof To check that the elements~(\ref{ttilde}) 
satisfy the defining relations (\ref{defrel}) 
we use Proposition~\ref{prop:qmr}. This proves the first part of the theorem.

To prove (\ref{qdettt}) we use induction on $m$.
For $m=1$ the identity is trivial. Let $m\geq 2$. We obtain from
(\ref{t:hatt}) that
\bea\label{T^(-1)}
\wh T(u)=\qdet T(u)\ts T(u-n+1)^{-1}.
\non
\eea
Taking the $mm$th entry on both sides we obtain
\bea\label{qdetmm}
\qdet  T(u)^{}_{{\cal C}{\cal C}}=\qdet T(u) \left(T(u-n+1)^{-1}\right)_{mm},
\eea
where ${\cal C}=\{1,\dots,m-1,m+1,\dots,n\}$. Similarly,
using the homomorphism (\ref{homtt}) we can write the
corresponding analog of (\ref{t:hatt}) for the matrix $\wt T(u)$
which gives
\bea\label{qdettildemm}
\qdet  \wt T(u)^{}_{{\cal A}'{\cal A}'}=\qdet \wt T(u)
\left(\wt T(u-m+1)^{-1}\right)_{mm}, 
\eea
where ${\cal A}'=\{1,\dots,m-1\}$. By the induction hypothesis,
\bea\label{qdetind}
\qdet  \wt T(u)^{}_{{\cal A}'{\cal A}'}=
\qdet  T(u)^{}_{{\cal C}{\cal C}} \ts
\qdet T(u-1)^{}_{{\cal B}{\cal B}}\cdots\qdet
T(u-m+2)^{}_{{\cal B}{\cal B}}. 
\non
\eea
It remains to apply (\ref{qdetmm}), (\ref{qdettildemm}) 
and the identity
\bea\label{qdetAA}
\left(T(u-n+1)^{-1}\right)^{}_{\cal A\cal A}=\wt T(u-m+1)^{-1}\ts 
\qdet T(u-m+1)^{}_{{\cal B}{\cal B}},
\non
\eea
which is a corollary of (\ref{t:hatt}) and (\ref{T:hatAA}).
\endproof

\section{Modified centralizer construction}	\label{mcc}
\setcounter{equation}{0}

We start with the definition of the Olshanski 
algebra $\A_m$; see~\cite{o:ea,o:ri}. 
Fix a nonnegative integer $m$ and
for any $n\geq m$ denote 
by $\g_m(n)$ the subalgebra in the general linear Lie algebra
$\gl(n)$ spanned by the basis elements $E_{ij}$ with 
$m+1\leq i,j\leq n$. The subalgebra $\g_m(n)$ is isomorphic to
$\gl(n-m)$.
Let $\A_m(n)$ denote the centralizer
of $\g_m(n)$ in the universal enveloping algebra $\A(n):=\U(\gl(n))$.
In particular, $\A_0(n)$ is the center of $\U(\gl(n))$.
Let $\A(n)^0$ denote the centralizer of $E_{nn}$
in $\A(n)$ and let $\I(n)$ be the left ideal in $\A(n)$ generated by the
elements $E_{in}$, $i=1,\dots,n$. Then $\I(n)^0:=\I(n)\cap \A(n)^0$
is a two-sided ideal in $\A(n)^0$ and one has a vector space
decomposition
\bea  \label{decomp}
\A(n)^0=\I(n)^0\oplus \A(n-1).
\non
\eea
Therefore the projection of $\A(n)^0$ onto $\A(n-1)$ 
with the kernel $\I(n)^0$ is an algebra homomorphism. Its
restriction to the subalgebra $\A_m(n)$ defines
a filtration preserving homomorphism
\bea\label{proj}
\pi_n: \A_m(n)\to \A_m(n-1)
\eea
so that one can define the algebra $\A_m$ as
the projective limit~(\ref{lim})
in the category of filtered algebras.
In other words, an element of the algebra $\A_m$ is a sequence
of the form
$a=(a_m,a_{m+1},\dots,a_n,\dots)$ where $a_n\in\A_m(n)$,
$\pi_{n}(a_n)=a_{n-1}$ for $n>m$, and
\bea  \label{deg}
\deg a:=\sup_{n\geq m}\deg a_n<\infty,
\non\eea
where
$\deg a_n$ denotes the degree of $a_n$ in  the universal
enveloping algebra $\A(n)$.

The algebra $\A_0$ is isomorphic to the algebra of virtual
Laplace operators for the Lie algebra $\gl(\infty)$~\cite{o:ri}.
Let $\lambda=(\lambda_1,\lambda_2,\dots)$ be an integer sequence
such that $\lambda_1\geq\lambda_2\geq\cdots$ and
$\lambda_i=0$ for all sufficiently large $i$ and let $L$ be
a representation of $\gl(\infty)$ with the highest weight
$\lambda$. Then every element $a\in \A_0$ acts 
in $L$ as a scalar. A family of algebraically independent generators
of $\A_0$ was constructed in~\cite[Remark~2.1.20]{o:ri}
(see also~\cite{gs:ci})  and their
eigenvalues in $L$ were calculated. We shall need
a different family of generators of $\A_0$ which will be more
suitable for our purposes.

Let $E=(E_{ij})$ denote the infinite matrix whose $ij$th entry is
the basis element $E_{ij}$ of $\gl(\infty)$ and let $E^{(n)}$
be its submatrix corresponding to the subset of indices
$1\leq i,j\leq n$. It is well-known
and can be easily verified that the mapping
\bea\label{homY(n)}
T(u)\mapsto 1+E^{(n)}u^{-1}
\eea
defines an algebra epimorphism
$\Y(n)\to\A(n)$.
Consider the quantum determinant
\bea  \label{qdet}
\qdet (1+E^{(n)}u^{-1})=\sum_{p\in\Sym_n}\sgn(p)\ts
(1+Eu^{-1})_{p(1)1}\cdots (1+E(u-n+1)^{-1})_{p(n)n}
\non\eea
and write
\bea  \label{coeff}
\qdet (1+E^{(n)}u^{-1})=1+{\E}_1^{(n)}u^{-1}+ {\E}_{2}^{(n)}u^{-2}+\cdots.
\non\eea
Since the coefficients of the series $\qdet T(u)$
are central in $\Y(n)$,
the elements ${\E}_i^{(n)}$ belong to the center
$\A_0(n)$ of $\A(n)$. 
Further, by definition~(\ref{proj}) we obviously have
\bea  \label{projqdet}
\pi_n: \qdet (1+E^{(n)}u^{-1})\mapsto \qdet (1+E^{(n-1)}u^{-1}).
\non\eea
So, we may define the elements ${\E}_i\in\A_0$, $i=1,2,\dots$ as
sequences ${\E}_i=({\E}_i^{(n)}|n\geq 1)$. We call the corresponding
generating series
\bea  \label{vqdet}
\qdet (1+Eu^{-1}):=1+{\E}_1u^{-1}+ {\E}_2u^{-2}+\cdots
\eea
the {\it virtual quantum determinant\/}. 
This is a particular case of
a more general construction of
quantum immanants given in \cite{o:qi,oo:ss}.
The following proposition
is immediate.

\bpr \label{prop:A_0}
The elements ${\E}_i$ with $i=1,2,\dots$ are algebraically independent
generators of the algebra $\A_0$. Their eigenvalues
in a representation $L$ of the Lie algebra $\gl(\infty)$
with the highest weight $\lambda$ are given by
the formula
\bea \label{eigen}
\qdet (1+Eu^{-1})|_L=\prod_{i=1}^{\infty}\frac{u+\lambda_i-i+1}{u-i+1}.
\non\eea
\epr

Consider the homomorphism $\Y(m)\to\Y(n)$ provided by Theorem~\ref{thm:qs}
and take its composition with (\ref{homY(n)}).
We obtain an algebra homomorphism $\Phi_n:\Y(m)\to \A(n)$
given by
\bea\label{Phi_n}
\Phi_n: t_{ij}(u)\mapsto \qdet (1+Eu^{-1})^{}_{{\cal B}_i{\cal B}_j},
\eea
where ${\cal B}_i$ denotes the set $\{i,m+1,\dots,n\}$.
In other words, $\Phi_n(t_{ij}(u))$ coincides with the image of the
element $\wt t_{ij}(u)$ under (\ref{homY(n)}). 
Proposition~\ref{prop:qmr} implies that this image commutes
with the elements of the subalgebra $\g_m(n)$ and so,
(\ref{Phi_n}) defines a homomorphism $\Phi_n:\Y(m)\to \A_m(n)$.
Furthermore, the family of homomorphisms $\{\Phi_n|n\geq m\}$
is obviously compatible with the homomorphisms (\ref{proj})
and thus defines an algebra homomorphism $\Phi: \Y(m)\to \A_m$.
Introducing the virtual quantum determinants of the
matrices $(1+Eu^{-1})^{}_{{\cal B}_i{\cal B}_j}$ with
${\cal B}_i=\{i,m+1,m+2,\dots\}$
as in
(\ref{vqdet}), we can represent this homomorphism as follows
\bea\label{Phi}
\Phi: t_{ij}(u)\mapsto \qdet (1+Eu^{-1})^{}_{{\cal B}_i{\cal B}_j}.
\non\eea

The following is a modified version of the Olshanski theorem;
see~\cite{o:ea, o:ri}. Denote by ${\wt \A}_0$ the commutative
algebra generated by the coefficients of the virtual quantum
determinant $\qdet (1+Eu^{-1})^{}_{{\cal B}{\cal B}}$ 
with ${\cal B}=\{m+1,m+2,\dots\}$.

\bth\label{thm:mcc}
The homomorphism $\Phi: \Y(m)\to \A_m$ is injective. Moreover, one has
an isomorphism
\bea\label{isom}
\A_m={\wt \A}_0\ot \Y(m),
\eea
where $\Y(m)$ is identified with its image under the embedding $\Phi$.
\eth

\Proof Consider the canonical filtration of the universal
enveloping algebra $\A(n)$. The corresponding
graded algebra $\gr \A(n)$ is isomorphic to
the symmetric algebra
$\P(n)$ of the space
$\gl(n)$. Elements of $\P(n)$ can be naturally identified
with polynomials in matrix elements of an $n\times n$-matrix
$x=(x_{ij})$.
Denote by
$\P_m(n)$ the subalgebra of the elements
of $\P(n)$ which are invariant under the adjoint action of
the Lie algebra $\g_m(n)$. One can show (see~\cite{o:ri}) that
the graded algebra $\gr \A_m$ is naturally isomorphic to
the projective limit $\P_m$ of the commutative algebras $\P_m(n)$
with respect to homomorphisms $\P_m(n)\to \P_m(n-1)$
analogous to (\ref{proj}). In particular, we can define the
virtual determinants
\bea\label{vdet} 
\det (1+x u^{-1})^{}_{{\cal B}_i{\cal B}_j},\qquad
\det (1+x u^{-1})^{}_{{\cal B}{\cal B}}
\eea
of the submatrices of the infinite
matrix $1+x u^{-1}$, where $x=(x_{ij})_{i,j=1}^{\infty}$,
as formal series in $u^{-1}$ with coefficients
from $\P_m$. These coefficients coincide with the images
in $\P_m$ of the coefficients of the virtual quantum determinants
of the corresponding submatrices of $E$. 
Therefore, both claims of the
theorem will follow from the fact that the 
coefficients of the series (\ref{vdet}) where $1\leq i,j\leq m$
are algebraically independent generators
of the algebra $\P_m$.
To see this, for a finite value of $n$ consider the 
sets ${\cal B}=\{m+1,m+2,\dots,n\}$ and ${\cal B}_i=\{i,m+1,m+2,\dots,n\}$.
We have an
identity
\bea\label{quasi}
\det (1+x u^{-1})^{}_{{\cal B}_i{\cal B}_j}
=\det (1+x u^{-1})^{}_{{\cal B}{\cal B}}\ts 
\left|(1+x u^{-1})^{}_{{\cal B}_i{\cal B}_j}\right|^{}_{ij}
\non\eea
with
\bea\label{paths}
\left|(1+x u^{-1})^{}_{{\cal B}_i{\cal B}_j}\right|^{}_{ij}=
\delta_{ij}+\sum_{k=1}^{\infty}(-1)^{k-1}\Lambda_{ij}^{(k)} u^{-k},
\non\eea
where 
$\Lambda_{ij}^{(k)}=\sum x_{ia_1}x_{a_1a_2}\cdots x_{a_{k-1}j}$, 
summed over the indices $a_r\in\{m+1,\dots,n\}$; 
see e.g.~\cite[Section~7.4]{gkllrt:ns}.
However, the coefficients of the polynomial 
$\det (1+x u^{-1})^{}_{{\cal B}{\cal B}}$ generate the full set
of invariants of the matrix $x^{}_{{\cal B}{\cal B}}$.
Further, these coefficients and the elements $\Lambda_{ij}^{(k)}$
generate the algebra $\P_m(n)$ which follows from
\cite[Section~2.1.10]{o:ri} and \cite[Section~7.4]{gkllrt:ns}.
Finally, given any positive integer $M$,
the coefficients at $u^{-1},\dots,u^{-M}$ of 
$\det (1+x u^{-1})^{}_{{\cal B}{\cal B}}$ and the
elements $\Lambda_{ij}^{(k)}$ with $k=1,\dots,M$
and $i,j=1,\dots,m$
are algebraically independent for
a sufficiently large $n$ which is 
implied by~\cite[Lemma~2.1.11]{o:ri}. \endproof


\noindent
{\bf Remark.} For $m> 1$ the algebra $\wt{\A}_0$ in 
the decomposition (\ref{isom})
can be replaced by $\A_0$. This can be deduced from the
classical Sylvester's identity (cf. Theorem~\ref{thm:qs}):
\bea\label{classy}
\det X=\det(1+x u^{-1})
\left(\det (1+x u^{-1})^{}_{{\cal B}{\cal B}}\right)^{m-1},
\non\eea
where $X=(X_{ij})$ with 
$X_{ij}=\det (1+x u^{-1})^{}_{{\cal B}_i{\cal B}_j}$. \endproof

\section{Extremal projection and transvector algebras}					
\label{tc}
\setcounter{equation}{0}

In the last two sections we shall consider the pair
of Lie algebras $\gl(m)\subset\gl(m+n)$ with the parameters
$m$ and $n$ fixed, where
$\gl(m)$ is spanned by the basis elements $E_{ij}$
of $\gl(m+n)$ with $i,j=1,\dots,m$.

Denote by $\h$ the diagonal Cartan subalgebra of
$\gl(m)$ spanned by $E_{ii}$ with
$i=1,\dots,m$. 
Consider the extension 
of the universal enveloping algebra $\A(m+n)=\U(\gl(m+n))$
\bea\label{ext}
\A'(m+n)=\A(m+n)\ot_{\U(\h)} \R(\h),
\non\eea
where $\R(\h)$ is the field of fractions of the commutative algebra
$\U(\h)$. Let $\J$
denote the left ideal in $\A'(m+n)$ generated by 
the elements $E_{ij}$ with $1\leq i<j\leq m$.
The {\it transvector algebra\/} $\Z=\Z(\gl(m+n),\gl(m))$ is the quotient algebra
of the normalizer
\bea\label{norm}
\Norm \J=\{x\in\A'(m+n)\ |\ \J x\subseteq \J\}
\non\eea
modulo the two-sided ideal $\J$; see~\cite{z:za}. 
It is an algebra over $\C$ and an $\R(\h)$-bimodule.
Generators of $\Z$
can be constructed by using the {\it extremal projection\/} $p$
for the Lie algebra $\gl(m)$ \cite{ast:po, z:sa}.
The projection $p$ is an element of an algebra $\F$
of formal series
and can be defined as follows. 
For any pair $(i,j)$
such that $1\leq i<j\leq m$ set
\bea\label{p_ab}
p_{ij}=\sum_{k=0}^{\infty}(E_{ji})^k(E_{ij})^k\ts
\frac{(-1)^k}{k!\ts (h_i-h_j+1)\cdots (h_i-h_j+k)},
\non\eea
where $h_i:=E_{ii}-i+1$. Then $p$ is given by
\bea\label{extrem}
p=\prod_{i<j}\ts p_{ij},
\non\eea
where the product is taken in any {\it normal ordering\/} on the pairs
$(i,j)$. The positive roots of 
$\gl(m)$ with respect to $\h$ (with the standard choice
of the positive root system)
are naturally enumerated by the pairs $(i,j)$
and normality means that
any composite root lies between its components.
The element $p$
is, up to a factor from $\R(\h)$,
a unique element
of $\F$ such that
\bea\label{annih}
E_{ij}\ts p=p E_{ji}=0\qquad{\rm for}\quad 1\leq i<j\leq m.
\eea
In particular, $p$ does not depend on the normal ordering
and
satisfies the conditions $p^2=p$ and $p^*=p$, where
$x\mapsto x^*$ is an involutive anti-automorphism of the algebra
$\F$ such that $(E_{ij})^*=E_{ji}$.

The action of $p$ on elements of the quotient $\A'(m+n)/\J$
is well-defined and the
transvector algebra $\Z$ 
can be identified with the image of $\A'(m+n)/\J$
with respect to the projection $p$:
\bea\label{p-image}
\Z=p\left(\A'(m+n)/\J\right).
\non\eea
An analog of the Poincar\'e--Birkhoff--Witt theorem holds for
the algebra $\Z$ \cite{z:za,z:it} so that ordered
monomials in the elements
$E_{ab},\ pE_{ia},\ pE_{ai}$,
where $i=1,\dots,m$ and $a,b=m+1,\dots,m+n$ form a basis of
$\Z$ as a left or right $\R(\h)$-module.
It will be convenient to use the following generators
of $\Z$: for $i=1,\dots,m$ and $a=m+1,\dots,m+n$
\bea\label{sima}
s_{ia}&=&pE_{ia}(h_i-h_1)\cdots (h_i-h_{i-1}),\non\\
\label{smai}
s_{ai}&=&pE_{ai}(h_i-h_{i+1})\cdots (h_i-h_m).
\eea
Explicitly,
\bea\label{simae}
s_{ia}&=&
\sum_{i>i_1>\cdots>i_s\geq 1}
E_{ii_1}E_{i_1i_2}\cdots E_{i_{s-1}i_s}E_{i_sa}
(h_i-h_{j_1})\cdots (h_i-h_{j_r}),\\
\label{s_maie}
s_{ai}&=&
\sum_{i<i_1<\cdots<i_s\leq m}
E_{i_1i}E_{i_2i_1}\cdots E_{i_{s}i_{s-1}}E_{ai_s}
(h_i-h_{j_1})\cdots (h_i-h_{j_r}),
\non
\eea
where $s=0,1,\dots$ and $\{j_1,\dots,j_r\}$ is the complementary subset
to $\{i_1,\dots,i_s\}$ respectively in the set $\{1,\dots,i-1\}$
or $\{i+1,\dots,m\}$. The elements $E_{ab}$, $a,b=m+1,\dots,m+n$ generate a subalgebra
in $\Z$ isomorphic to $\A(n)$. Moreover, one also has
the following relations in $\Z$
which are obtained by using the methods of~\cite{z:za,z:it}.
Let the indices $i,j$ and $a,b,c$ run through
the sets $\{1,\dots,m\}$ and $\{m+1,\dots,m+n\}$, respectively. We have
\bea\label{Ess}
[E_{ab}, s_{ci}]=\delta_{bc}s_{ai},\qquad
[E_{ab}, s_{ic}]=-\delta_{ac}s_{ib}.
\eea
Moreover, if $i\ne j$ then
\bea\label{s_ias_jb}
s_{ai}s_{bj}=s_{bj}s_{ai}\ts\frac{h_i-h_j+1}{h_i-h_j}-
s_{bi}s_{aj}\ts\frac{1}{h_i-h_j},
\non\eea
while
\bea\label{s_ias_ja}
s_{ai}s_{bi}=s_{bi}s_{ai}.
\eea
Finally, $s_{ia}s_{bj}=s_{bj}s_{ia}$ for $i\ne j$ while
\beq\label{s_ias_bi}
s_{ia}s_{bi}=(\delta_{ab}(E_{ii}+m-i)-E_{ab})\prod_{j=1,\ts j\ne i}^m(h_i-h_j-1)
+\sum_{j=1}^m s_{bj}s_{ja}
\prod_{k=1,\ts k\ne j}^m \frac{h_i-h_k-1}{h_j-h_k}.
\eeq

The anti-involution $x\mapsto x^*$ extends to the algebra $\Z$
so that $(pE_{ai})^*=pE_{ia}$~\cite{z:za,z:it}. Therefore,
\beq\label{conj}
(s_{ai})^*=s_{ia}\frac{(h_i-h_{i+1}+1)\cdots (h_i-h_{m}+1)}
{(h_i-h_1)\cdots (h_i-h_{i-1})}.
\eeq

The centralizer $\A=\U(\gl(m+n))^{\gl(m)}$ is a subalgebra in the normalizer
$\Norm \J$ and so, one has the natural algebra homomorphism
$\nu:\A\to \Z$. 
By the results of Section~\ref{mcc} we have an algebra
homomorphism
\beq\label{homYA}
\phi:\Y(n)\to \A,\qquad
t_{ab}(u)\mapsto \qdet (1+Eu^{-1})^{}_{{\cal C}_a{\cal C}_b},
\eeq
where ${\cal C}_a=\{1,\dots,m,m+a\}$ for $a=1,\dots,n$; see~(\ref{Phi_n}).
The composition $\psi=\nu\circ\phi$
is an algebra homomorphism $\psi:\Y(n)\to\Z$. The next theorem provides
explicit formulas for the images of the generators of $\Y(n)$
with respect to $\psi$. 
Note that by (\ref{homYA}) the image of
\bea\label{T_ij}
T_{m+a,m+b}(u):=u(u-1)\cdots (u-m)\ts t_{ab}(u)
\eea
under $\psi$ is a polynomial in $u$.

\bth\label{thm:images}
We have for $a,b=m+1,\dots,m+n$:
\beq\label{psiTab}
\psi:T_{ab}(u)\mapsto (\delta_{ab}(u-m)+E_{ab})
\prod_{i=1}^m(u+h_i)
-\sum_{i=1}^m s_{ai}s_{ib}
\prod_{j=1,\ts j\ne i}^m \frac{u+h_j}{h_i-h_j}
\eeq
and
\beq\label{psiTab2}
\psi:T_{ab}(u)\mapsto (\delta_{ab}u+E_{ab})
\prod_{i=1}^m(u+h_i-1)
-\sum_{i=1}^m s_{ib}s_{ai}
\prod_{j=1,\ts j\ne i}^m \frac{u+h_j-1}{h_i-h_j}.
\eeq
\eth

\Proof By (\ref{qminor})
the image of $T_{ab}(u)$ under the homomorphism (\ref{homYA}) is
\beq\label{qdetpol}
\sum_{\sigma}\sgn(\sigma)\ts
(u+E)_{\sigma(1),1}\cdots(u+E-m+1)_{\sigma(m),m}(u+E-m)_{\sigma(a),b},
\eeq
summed over permutations $\sigma$ of the set $\{1,\dots,m,a\}$.
To find the image of the polynomial (\ref{qdetpol}) in $\Z$ we shall
regard it modulo the ideal $\J$ and apply the projection $p$.
Consider first the sum in (\ref{qdetpol}) 
over the permutations $\sigma$ such that $\sigma(a)=a$. Using
(\ref{annih}) we obtain that its image
in $\Z[u]$ is
\bea\label{highest}
(u+h_1)\cdots (u+h_m)(\delta_{ab}(u-m)+E_{ab}).
\eea
Further, consider the remaining summands in (\ref{qdetpol}) 
and suppose that $\sigma(a)=i$, $\sigma(j)=a$ for certain
$i,j\in\{1,\dots,m\}$. Using the property (\ref{annih}) of 
the projection $p$ we find that the image of such a summand 
in $\Z[u]$ can be
nonzero only if $i\geq j$ and $\sigma$ is a cyclic
permutation of the form $\sigma=(a,i_1,\dots,i_s)$ where
$i=i_1>i_2>\cdots>i_s=j$. We have the equality modulo $\J$:
\bea\label{p-act}
E_{ai_s}E_{i_si_{s-1}}\cdots E_{i_2i_1}E_{i_1b}
=E_{aj}E_{jb}.
\non\eea
Since $\sgn(\sigma)=(-1)^s$
the image in $\Z[u]$ of the summand corresponding to
$\sigma$ is
\bea\label{imsigma}
(-1)^s\ts pE_{aj}E_{jb}(u+h_1)\cdots \wh{(u+h_{i_s})}\cdots
\wh{(u+h_{i_1})}\cdots (u+h_m),
\non\eea
where the hats indicate the factors to be omitted. Taking the sum over
cycles $\sigma$ we find that the polynomial $S_{ab}(u):=\psi(T_{ab}(u))$
is given by
\bea\label{Z_ijprem}
S_{ab}(u)&=&(\delta_{ab}(u-m)+E_{ab})(u+h_1)\cdots (u+h_m)\non\\
\mbox{}&-&\sum_{j=1}^m 
pE_{aj}E_{jb} (u+h_1)\cdots (u+h_{j-1})
(u+h_{j+1}-1)\cdots (u+h_m-1). \non
\eea
Consider the value of $S_{ab}(u)$ at $u=-h_i$ for 
$i\in\{1,\dots,m\}$ (to make this evaluation well-defined we agree to write
the coefficients of the polynomial to the left of the powers of $u$).
We obtain
\[
S_{ab}(-h_i)=-\sum_{j=1}^i 
pE_{aj}E_{jb} (h_1-h_i)\cdots (h_{j-1}-h_i)
(h_{j+1}-h_i-1)\cdots (h_m-h_i-1).
\]
On the other hand, by 
(\ref{smai}) we have
\[
s_{ai}s_{ib}=pE_{ai}\ts s_{ib}
(h_i-h_{i+1}+1)\cdots (h_i-h_m+1).
\]
Finally, using
(\ref{simae}) and (\ref{annih})
we verify the identity 
$s_{ai}s_{ib}=(-1)^{m}S_{ab}(-h_i)$.
The difference of $S_{ab}(u)$ and the product (\ref{highest})
is a polynomial in $u$ of degree $< m$. So, (\ref{psiTab})
follows from the Lagrange interpolation formula. 

The relation (\ref{psiTab2}) follows from (\ref{s_ias_bi})
and (\ref{psiTab}); it suffices to consider the values
of the polynomial $S_{ab}(u)$ at $u=-h_i+1$ for $i=1,\dots,m$.
It can also be proved by the above argument
where (\ref{qminor}) is applied to
${t\ts}^{a,1\cdots\ts m}_{b,1\cdots\ts m}(u)$
(thus giving another derivation of (\ref{s_ias_bi})).
\endproof


\section{Elementary representations of the Yangian}    \label{ery}
\setcounter{equation}{0}

Here we use the results of the previous section
to study a special class of representations
of the Yangian $\Y(n)$ called elementary.
They play an important role
in the classification of the representations of $\Y(n)$
with a semisimple action of 
the Gelfand--Tsetlin subalgebra; see~\cite{nt:ry}.

\subsection{Yangian action on the multiplicity space}

A representation $L$
of the Yangian $\Y(n)$ 
is called {\it highest weight\/}
if it is generated by a nonzero vector $\zeta$ (the {\it highest vector\/})
such that
\bea\label{hvY}
t_{aa}(u)\zeta=\lambda_a&
\displaystyle{(u)\zeta\qquad}&{\rm for}\quad a=1,\dots,n,\non\\
t_{ab}(u)\zeta=0\ts&
\displaystyle{\qquad}&{\rm for}\quad 1\leq a<b\leq n
\non\eea
for certain formal series $\lambda_a(u)\in 1+u^{-1}\C[[u^{-1}]]$.
The set $\lambda(u)=(\lambda_1(u),\dots,\lambda_n(u))$ is called
the {\it highest weight\/} of $L$; cf. \cite{cp:yr,d:nr}. 
Every finite-dimensional irreducible representation
of the Yangian $\Y(n)$ is highest weight. It contains a unique,
up to scalar multiples, highest vector.
An irreducible representation of $\Y(n)$ with
the highest weight $\lambda(u)$
is finite-dimensional
if and only if there exist monic polynomials
$P_1(u),\dots,P_{n-1}(u)$ in $u$ 
(called the {\it Drinfeld polynomials\/}) such that
\bea\label{drpol}
\frac{\lambda_a(u)}{\lambda_{a+1}(u)}=
\frac{P_a(u+1)}{P_a(u)},\qquad a=1,\dots,n-1.
\eea
These results are contained in \cite{d:nr}; see also \cite{cp:yr,m:fd}.

Let $\lambda=(\lambda_1,\dots,\lambda_{m+n})$ be an
$(m+n)$-tuple of complex numbers satisfying the condition
$\lambda_i-\lambda_{i+1}\in\ZZ_+$
for $i=1,\dots,m+n-1$. Denote by
$L(\lambda)$ the irreducible finite-dimensional representation of the
Lie algebra $\gl(m+n)$ with the highest weight $\lambda$.
It contains a unique
nonzero vector $\xi$ (the highest vector) such that
\bea\label{hv}
E_{ii}\xi=\lambda_i&
\displaystyle{\xi\qquad}&{\rm for}\quad i=1,\dots,m+n,\non\\
E_{ij}\xi=0\ts&\displaystyle{\qquad}&{\rm for}\quad 1\leq i<j\leq m+n.
\non\eea

Given a $\gl(m)$-highest weight $\mu=(\mu_1,\dots,\mu_m)$
we denote by $L(\lambda)^+_{\mu}$ the subspace of
$\gl(m)$-highest vectors in $L(\lambda)$ of weight $\mu$:
\bea\label{hvglm}
L(\lambda)^+_{\mu}=\{\eta\in L(\lambda)\ts|\ts E_{ii}\eta=
\mu_i&
\displaystyle{\eta
\qquad}&{\rm for}\quad i=1,\dots,m,\non\\
E_{ij}\eta=0\ts&
\displaystyle{\qquad}&{\rm for}\quad 1\leq i<j\leq m\}.
\non\eea
The dimension of $L(\lambda)^+_{\mu}$ coincides with
the multiplicity of the $\gl(m)$-module
$L(\mu)$ in 
the restriction of
$L(\lambda)$ to $\gl(m)$.
The multiplicity
space $L(\lambda)^+_{\mu}$ admits a natural structure of
an irreducible representation of the centralizer algebra
$\A$~\cite[Section~9.1]{d:ae}. On the other hand, the homomorphism
(\ref{homYA}) allows us to regard $L(\lambda)^+_{\mu}$
as a $\Y(n)$-module.
As it follows from the proof of Theorem~\ref{thm:mcc},
the algebra $\A$ is generated by the image of $\phi$ and
the center of $\U(\gl(m))$. 
The elements of the center act in
$L(\lambda)^+_{\mu}$ as scalar operators and so, the $\Y(n)$-module
$L(\lambda)^+_{\mu}$
is irreducible. Following~\cite{nt:ry} we call it
{\it elementary}.

We shall use a {\it contravariant bilinear form\/} $\langle\ts ,\rangle$
on the space $L(\lambda)$ which is uniquely determined by
the conditions
\bea\label{form} 
\langle \xi,\xi\rangle=1;\qquad\langle Xu,v\rangle =\langle u,X^*v\rangle
\quad {\rm for}\quad u,v\in L(\lambda),\quad X\in\A(m+n),
\non\eea
where $X\mapsto X^*$ is the anti-involution on $\A(m+n)$
such that $(E_{ab})^*=E_{ba}$.
The form $\langle\ts ,\rangle$ is nondegenerate on $L(\lambda)$,
and so is its restriction to each nonzero subspace $L(\lambda)^+_{\mu}$. 
Indeed, if there is a vector $\eta\in L(\lambda)^+_{\mu}$ such that
$\langle \eta,\eta'\rangle=0$ for each $\eta'\in L(\lambda)^+_{\mu}$
then $\eta$ is also orthogonal to all elements of $L(\lambda)^+$
because distinct $\gl(m)$-weight subspaces in $L(\lambda)^+$
are pairwise orthogonal to each other. On the other hand,
$L(\lambda)=\U({\mathfrak n}_-)\ts L(\lambda)^+$, where 
${\mathfrak n}_-$ is the lower triangular subalgebra in $\gl(m)$.
Therefore, $\eta$ is orthogonal to $L(\lambda)$ since for any 
$X\in {\mathfrak n}_-$ and $\eta'\in L(\lambda)$ one has
\bea\label{eta}
\langle \eta,X\eta'\rangle=\langle X^*\eta,\eta'\rangle=0
\non\eea
because $X^*\eta=0$. Thus, $\eta=0$.

We shall maintain the notation $\langle\ts ,\rangle$ for the restriction
of the form to a subspace $L(\lambda)^+_{\mu}$. 

\subsection{Quantum minor formulas for generators of $\Z$}

Using the homomorphism $\Y(m+n)\to \A(m+n)$ defined by
(\ref{homY(n)}) we can carry over the
quantum minor identities from Section~\ref{qst}
to the matrix $1+Eu^{-1}$ or,
which is more convenient, to
$E(u)=u+E$.
In accordance with (\ref{qminor}) and (\ref{qminor2}) we
define the polynomials 
${E\ts}^{a_1\cdots\ts a_s}_{b_1\cdots\ts b_s}(u)$
by the following equivalent formulas
\bea\label{qminorE}
{E\ts}^{a_1\cdots\ts a_s}_{b_1\cdots\ts b_s}(u)&=&
\sum_{\sigma\in \Sym_s} \sgn(\sigma)\ts (u+E)_{a_{\sigma(1)}b_1}\cdots
(u+E-s+1)_{a_{\sigma(s)}b_s}\\
\label{qminorE2}
\mbox{}&=&
\sum_{\sigma\in \Sym_s} \sgn(\sigma)\ts (u+E-s+1)_{a_1b_{\sigma(1)}}\cdots
(u+E)_{a_sb_{\sigma(s)}}.
\eea
The polynomial ${E\ts}^{a_1\cdots\ts a_s}_{b_1\cdots\ts b_s}(u)$
is skew-symmetric with respect to
permutations of the upper indices and of the lower indices.
We shall need to evaluate $u$ in $\R(\h)$. To make this 
operation well-defined let us agree to write the coefficients
of this polynomial to the left from powers of $u$.
In the following proposition we use the notation from Section~\ref{tc}.

\bpr\label{prop:eq}
For any $i\in\{1,\dots,m\}$ and $a\in\{m+1,\dots,m+n\}$ we have
the equalities in $\A'(m+n)$ modulo the ideal $\J$:
\bea\label{siaE}
s_{ia}&=&(-1)^{i-1}{E\ts}^{1\cdots\ts i}_{1\cdots\ts i-1,a}(-h_i),\\
\label{saiE}
s_{ai}&=&{E\ts}^{i+1\cdots\ts m,a}_{i\cdots\ts m}(-h_i-i+1).
\eea
\epr

\Proof Let us show first that 
${E\ts}^{i+1\cdots\ts m,a}_{i\cdots\ts m}(-E_{ii})$ 
(recall that $h_i=E_{ii}-i+1$) belongs
to the normalizer $\Norm \J$ of $\J$.
By Proposition~\ref{prop:qmr} we have the relations
\bea\label{EE}
[E_{ij}, {E\ts}^{a_1\cdots\ts a_s}_{b_1\cdots\ts b_s}(u)]=
\sum_{r=1}^s\left(\delta_{ja_r}
{E\ts}^{a_1\cdots\ts i\cdots\ts a_s}_{b_1\cdots\ts b_s}(u)
-\delta_{ib_r}
{E\ts}^{a_1\cdots\ts a_s}_{b_1\cdots\ts j\cdots\ts b_s}(u)
\right),
\non\eea
where $i$ and $j$ on the right hand side take the $r$th positions.
This implies that $E_{k,k+1}$ commutes with
${E\ts}^{i+1\cdots\ts m,a}_{i\cdots\ts m}(u)$ for 
$k\in\{1,\dots,m-1\}$,
$k\ne i$. Hence,
\bea\label{EkE}
E_{k,k+1}{E\ts}^{i+1\cdots\ts m,a}_{i\cdots\ts m}(-E_{ii})\in\J.
\non\eea
Furthermore,
\bea\label{EEcom}
[E_{i,i+1}, {E\ts}^{i+1\cdots\ts m,a}_{i\cdots\ts m}(u)]=
{E\ts}^{i,i+2\cdots\ts m,a}_{i\cdots\ts m}(u).
\non\eea
We have 
\bea\label{EpmE}
{E\ts}^{i,i+2\cdots\ts m,a}_{i\cdots\ts m}(u)=
(-1)^{m-i}{E\ts}^{i+2\cdots\ts m,a,i}_{i\cdots\ts m}(u),
\non\eea
and by (\ref{qminorE2})
\bea\label{qmE}
{E\ts}^{i+2\cdots\ts m,a,i}_{i\cdots\ts m}(u)=
\sum_{\sigma} \sgn(\sigma)\ts (u+E-m+i)_{i+2,\sigma(i)}\cdots
(u+E)_{i,\sigma(m)}.
\eea
However, $\sigma(m)$ take values in
$\{i,\dots,m\}$
and so, the element (\ref{qmE})
belongs to the ideal $\J$
for $u=-E_{ii}$.
Note that $E_{ii}$ belongs to the normalizer $\Norm \J$
and therefore
${E\ts}^{i+1\cdots\ts m,a}_{i\cdots\ts m}(-E_{ii})$
does.

To complete the proof of (\ref{saiE}) we calculate the image
${E\ts}^{i+1\cdots\ts m,a}_{i\cdots\ts m}(-E_{ii})$
under the extremal projection $p$. By (\ref{qminorE})
we have
\bea\label{qmEE}
{E\ts}^{i+1\cdots\ts m,a}_{i\cdots\ts m}(u)=
\sum_{\sigma} \sgn(\sigma)\ts (u+E)_{\sigma(i+1),i}\cdots
(u+E-m+i)_{\sigma(a),m}.
\eea
By the property (\ref{annih}) of $p$ we obtain that
the image of a summand in (\ref{qmEE}) under $p$
is zero unless
$\sigma(i+1)=a$. Further, a summand is zero modulo the ideal $\J$
unless
\bea\label{sigma} 
\sigma(a)=m,\quad\sigma(m)=m-1,\quad\ldots\ts,\quad \sigma(i+2)=i+1.
\non\eea
Therefore,
the image of (\ref{qmEE}) with $u=-E_{ii}$ equals
\bea\label{pEE}
p E_{ai}\ts(h_i-h_{i+1})\cdots (h_i-h_m),
\non\eea
which coincides with $s_{ai}$. This proves (\ref{saiE}).
The relation (\ref{siaE}) is proved by a similar argument. \endproof

\subsection{Highest vector of the $\Y(n)$-module $L(\lambda)^+_{\mu}$}

Our aim now is to explicitly construct the highest vector
of the $\Y(n)$-module $L(\lambda)^+_{\mu}$ and to find
its highest weight. This will also allow us to calculate its
Drinfeld polynomials.

The action of the Yangian $\Y(n)$ 
in $L(\lambda)^+_{\mu}$ is determined by
Theorem~\ref{thm:images}. From now on we shall assume
that the highest weight $\lambda$ is a partition, that is,
the $\lambda_i$ are nonnegative integers. This does not lead
to a real loss of generality because the formulas and
arguments below can be easily adjusted to 
be valid in the general case. Given a general $\lambda$ one
can add a suitable complex number to all entries of $\lambda$
to get a partition.

As it follows from the branching rule for 
the general linear Lie algebras (see \cite{z:cl})
the space $L(\lambda)^+_{\mu}$ is nonzero only if $\mu$
is a partition such that $\mu\subset\lambda$ and each column of the skew
diagram $\lambda/\mu$ contains at most $n$ cells. The figure
below illustrates the skew diagram for $\lambda=(10,8,5,4,2)$
and $\mu=(6,3)$.

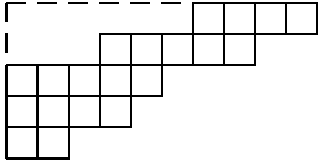
\setlength{\unitlength}{0.5em}
\begin{center}
\begin{picture}(20,10)

\put(12,10){\line(1,0){8}}
\put(6,8){\line(1,0){14}}
\put(0,6){\line(1,0){16}}
\put(10,6){\line(1,0){6}}
\put(0,4){\line(1,0){10}}
\put(0,2){\line(1,0){8}}
\put(0,0){\line(1,0){4}}

\put(0,0){\line(0,1){6}}
\put(2,0){\line(0,1){6}}
\put(4,0){\line(0,1){6}}
\put(6,2){\line(0,1){6}}
\put(8,2){\line(0,1){6}}
\put(10,4){\line(0,1){4}}
\put(12,6){\line(0,1){4}}
\put(14,6){\line(0,1){4}}
\put(16,6){\line(0,1){4}}
\put(18,8){\line(0,1){2}}
\put(20,8){\line(0,1){2}}

\multiput(0,10)(2,0){6}{\line(1,0){1.2}}

\multiput(0,10)(0,-2){2}{\line(0,-1){1.2}}

\end{picture}
\end{center}
\setlength{\unitlength}{1pt}

\noindent
Introduce the {\it row order\/} 
on the cells of $\lambda/\mu$ corresponding
to reading the diagram by rows from left to right starting from
the top row. For a cell $\alpha\in\lambda/\mu$ denote
by $r(\alpha)$ the row number of $\alpha$ and by $l(\alpha)$
the (increased) leglength of $\alpha$ which equals 1 plus
the number of cells of $\lambda/\mu$ in the column containing $\alpha$
which are below $\alpha$. Consider the following element
of $L(\lambda)$:
\bea\label{zetas}
\zeta=\prod_{\alpha\in\lambda/\mu,\  r(\alpha)\leq m} 
s^{}_{m+l(\alpha),r(\alpha)}\ts  \xi,
\eea
where $\xi$ is the highest vector of $L(\lambda)$ and the product is
taken in the row order. Using (\ref{Ess}) we find that
$\zeta\in L(\lambda)^+_{\mu}$.
For the above example of $\lambda/\mu$ we have
$m=2$, $n=3$ and
\bea\label{zetaex}
\zeta=(s^{}_{41})^2\ts (s^{}_{31})^2\ts 
s^{}_{52}\ts s^{}_{42}\ts (s^{}_{32})^3\ts\xi.
\non\eea
Proposition~\ref{prop:eq} allows us to rewrite
the definition (\ref{zetas})
in a different form. Introduce the notation
\bea\label{tau}
\tau_{ai}(u)={E\ts}^{i+1\cdots\ts m,a}_{i\cdots\ts m}(u),
\qquad
\tau_{ia}(u)={E\ts}^{1\cdots\ts i}_{1\cdots\ts i-1,a}(u).
\non\eea
Then by (\ref{saiE}) we have
\bea\label{zetatau}
\zeta=\prod_{\alpha\in\lambda/\mu,\  r(\alpha)\leq m} 
\tau^{}_{m+l(\alpha),r(\alpha)}(-c(\alpha))\ts \xi,
\non\eea
where $c(\alpha)$ is the column number of $\alpha$, and the product
is taken in the row order. 

Given three integers $i,j,k$
we shall denote by $\middle\{i,j,k\}$ that of the three which is
between the two others.

\bth\label{thm:hvpr}
The vector $\zeta$ defined by 
(\ref{zetas})
is the highest vector of the $\Y(n)$-module $L(\lambda)^+_{\mu}$.
The highest weight of this module is
$(\lambda_{1}(u),\dots,\lambda_{n}(u))$ where
\beq\label{lambda}
\lambda_a(u)=\frac{(u+\nu^{(1)}_a)(u+\nu^{(2)}_a-1)\cdots
(u+\nu^{(m+1)}_a-m)}{u(u-1)\cdots (u-m)}
\eeq
and
\beq\label{nu}
\nu^{(i)}_a=\middle\{\mu_{i-1},\mu_i,\lambda_{a+i-1}\}
\eeq
with $\mu_{m+1}=0$, and $\mu_0$ is considered to be sufficiently
large.
\eth

Note that for each $i$
the $n$-tuple $\nu^{(i)}=(\nu^{(i)}_{1},\dots,\nu^{(i)}_{n})$
is a partition which can be obtained from $\lambda/\mu$
as follows. Consider the subdiagram
of $\lambda$ of the form
$\lambda^{(i)}=(\lambda_i,\lambda_{i+1},\dots,\lambda_{i+n-1})$. Replace
the rows of $\lambda^{(i)}$ which are longer than $\mu_{i-1}$ 
by $\mu_{i-1}$ while those which are shorter than $\mu_i$ replace with
$\mu_i$ and leave 
the remaining rows unchanged. 
The resulting partition is $\nu^{(i)}$.
For the above example
with $\lambda=(10,8,5,4,2)$
and $\mu=(6,3)$ we have
\bea\label{nu123}
\nu^{(1)}=(10,8,6),
\qquad
\nu^{(2)}=(6,5,4),
\qquad
\nu^{(3)}=(3,3,2),
\non\eea
as illustrated:

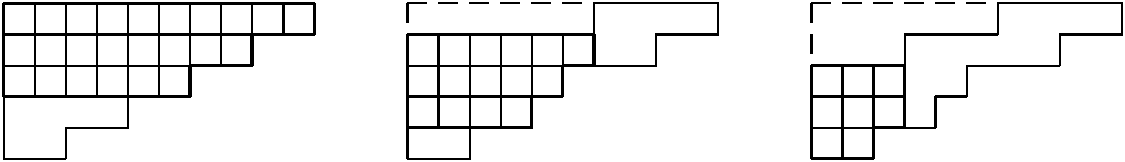
\setlength{\unitlength}{0.5em}
\begin{center}
\begin{picture}(72,10)

\put(12,10){\line(1,0){8}}
\put(0,8){\line(1,0){16}}
\put(0,6){\line(1,0){12}}
\put(10,6){\line(1,0){6}}
\put(4,2){\line(1,0){4}}
\put(0,0){\line(1,0){4}}

\put(0,0){\line(0,1){4}}
\put(2,4){\line(0,1){6}}
\put(4,0){\line(0,1){2}}
\put(4,4){\line(0,1){6}}
\put(6,4){\line(0,1){6}}
\put(8,2){\line(0,1){8}}
\put(10,4){\line(0,1){6}}
\put(12,6){\line(0,1){4}}
\put(14,6){\line(0,1){4}}
\put(16,6){\line(0,1){4}}
\put(18,8){\line(0,1){2}}
\put(20,8){\line(0,1){2}}

\multiput(0,10)(2,0){6}{\line(1,0){1.2}}

\multiput(0,10)(0,-2){2}{\line(0,-1){1.2}}

\thicklines

\put(0,10){\line(1,0){20}}
\put(16,8){\line(1,0){4}}
\put(12,6){\line(1,0){4}}
\put(0,4){\line(1,0){12}}

\put(0,4){\line(0,1){6}}
\put(12,4){\line(0,1){2}}
\put(16,6){\line(0,1){2}}
\put(20,8){\line(0,1){2}}

\thinlines

\put(26,0){\begin{picture}(20,10)

\put(12,10){\line(1,0){8}}
\put(16,8){\line(1,0){4}}
\put(0,6){\line(1,0){16}}
\put(10,6){\line(1,0){6}}
\put(0,4){\line(1,0){10}}
\put(0,2){\line(1,0){8}}
\put(0,0){\line(1,0){4}}

\put(0,0){\line(0,1){6}}
\put(2,2){\line(0,1){6}}
\put(4,0){\line(0,1){8}}
\put(6,2){\line(0,1){6}}
\put(8,2){\line(0,1){6}}
\put(10,4){\line(0,1){4}}
\put(12,6){\line(0,1){4}}
\put(16,6){\line(0,1){2}}
\put(20,8){\line(0,1){2}}

\multiput(0,10)(2,0){6}{\line(1,0){1.2}}

\multiput(0,10)(0,-2){2}{\line(0,-1){1.2}}

\thicklines

\put(0,8){\line(1,0){12}}
\put(10,6){\line(1,0){2}}
\put(8,4){\line(1,0){2}}
\put(0,2){\line(1,0){8}}

\put(0,2){\line(0,1){6}}
\put(8,2){\line(0,1){2}}
\put(10,4){\line(0,1){2}}
\put(12,6){\line(0,1){2}}

\thinlines

\end{picture}
}

\put(52,0){\begin{picture}(20,10)

\put(12,10){\line(1,0){8}}
\put(6,8){\line(1,0){6}}
\put(16,8){\line(1,0){4}}
\put(10,6){\line(1,0){6}}
\put(8,4){\line(1,0){2}}
\put(0,4){\line(1,0){6}}
\put(0,2){\line(1,0){8}}
\put(0,0){\line(1,0){4}}

\put(0,0){\line(0,1){6}}
\put(2,0){\line(0,1){6}}
\put(4,0){\line(0,1){6}}
\put(6,2){\line(0,1){6}}
\put(8,2){\line(0,1){2}}
\put(10,4){\line(0,1){2}}
\put(12,8){\line(0,1){2}}
\put(16,6){\line(0,1){2}}
\put(20,8){\line(0,1){2}}

\multiput(0,10)(2,0){6}{\line(1,0){1.2}}

\multiput(0,10)(0,-2){2}{\line(0,-1){1.2}}

\thicklines

\put(0,6){\line(1,0){6}}
\put(4,2){\line(1,0){2}}
\put(0,0){\line(1,0){4}}

\put(0,0){\line(0,1){6}}
\put(4,0){\line(0,1){2}}
\put(6,2){\line(0,1){4}}

\thinlines

\end{picture}
}

\end{picture}
\end{center}
\setlength{\unitlength}{1pt}

\noindent{\bf Proof of Theorem \ref{thm:hvpr}.\ \ }
Introduce two parameters $k$ and $l$ of the diagram $\lambda/\mu$
as follows. We let $k$ be the row number of the top (non-empty)
row of $\lambda/\mu$ if this number is less or equal to $m$;
otherwise set $k=m+1$. So, if $k\leq m$ then we have
$\mu_i=\lambda_i$ for 
$i=1,\dots,k-1$ and $\mu_k<\lambda_k$. 
We denote by $l$
the leglength of the cell $\alpha=(k,\mu_k+1)$
with $k\leq m$.
We shall prove 
the following three relations simultaneously by
induction on $k$ and $l$
(see (\ref{T_ij})
for the definition of $T_{ab}(u)$):
\bea\label{Tabzeta}
T_{ab}(u)\ts\zeta=0
\eea
for $m+1\leq a<b\leq m+n$,
\beq\label{Taazeta}
T_{m+a,m+a}(u)\ts\zeta=(u+\nu^{(1)}_a)(u+\nu^{(2)}_a-1)\cdots
(u+\nu^{(m+1)}_a-m)\ts\zeta
\eeq
for $a=1,\dots, n$, and
\beq\label{szeta}
s_{m+l,k}\ts\zeta=0,
\eeq
if $k\leq m$ and $\mu_k=\lambda_{k+l}$.

Suppose that $k=m+1$. Then $\zeta=\xi$ and both 
(\ref{Tabzeta}) and (\ref{Taazeta}) are immediate from 
(\ref{psiTab}). Now let $k\leq m$ be fixed. 
We assume that (\ref{Tabzeta})--(\ref{szeta}) hold for
all greater values of the parameter $k$.
We proceed by induction on $l$.
 
To simplify the notation
we shall identify the series $T_{ab}(u)$
with its image under the homomorphism (\ref{homYA}). That is, we set
\bea\label{TabE}
T_{ab}(u)={E\ts}^{1\cdots\ts m,a}_{1\cdots\ts m,b}(u).
\non\eea
We have
\beq\label{T*}
T_{ab}(u)^*=T_{ba}(u),
\eeq
which follows from (\ref{qminorE}) and
(\ref{qminorE2}).
The following relations are derived from Proposition~\ref{prop:qmr}:
\bea\label{Ttau}
T_{ab}(u)\tau_{ci}(v)&=&\tau_{ci}(v)T_{ab}(u)\ts
\frac{u-v-i}{u-v-i+1}+\tau_{ai}(v)T_{cb}(u)\ts
\frac{1}{u-v-i+1},\\
\label{Ttau2}
T_{ab}(u)\tau_{ci}(v)&=&\tau_{ci}(v)T_{ab}(u)\ts
\frac{u-v-i-1}{u-v-i}+T_{cb}(u)\tau_{ai}(v)\ts
\frac{1}{u-v-i}.
\eea

Using (\ref{saiE}) replace the first factor $s_{m+l,k}$ in
(\ref{zetas}) by $\tau_{m+l,k}(-\mu_k-1)$ so that the vector $\zeta$
will be written in the form $\zeta=\tau_{m+l,k}(-\mu_k-1)\zeta'$. 
Apply the operator $T_{ab}(u)$ with $a<b$ to $\zeta$. If $b>m+l$
then by (\ref{Ttau}) and (\ref{saiE}) we have
\beq\label{Ttaumu}
T_{ab}(u)\zeta=
\frac{u+\mu_k-k+1}{u+\mu_k-k+2}\ts s_{m+l,k}T_{ab}(u)\ts\zeta'
+\frac{1}{u+\mu_k-k+2}\ts s_{ak}
T_{m+l,b}(u)\ts\zeta'.
\eeq
Now (\ref{Tabzeta}) follows by induction on $l$ and the degree
of $\zeta$ with respect to $s_{m+l,k}$.

If $b\leq m+l$ then $a<m+l$ and by (\ref{Ttau2}) and (\ref{saiE})
we have
\beq\label{Ttaumu2}
T_{ab}(u)\zeta=
\frac{u+\mu_k-k}{u+\mu_k-k+1}\ts s_{m+l,k}T_{ab}(u)\ts\zeta'
+\frac{1}{u+\mu_k-k+1}\ts T_{m+l,b}(u)s_{ak}\ts\zeta'.
\eeq
However, by (\ref{s_ias_ja})
we have 
$s_{ak}\ts\zeta'=0$, where we have used
the induction hypothesis for (\ref{szeta}).
So, 
(\ref{Tabzeta}) follows again by an obvious induction.

Now we use a similar argument to calculate
$T_{aa}(u)\ts\zeta$. The relation (\ref{Ttaumu}) with $a=b$
gives
\bea\label{Taamu}
T_{aa}(u)\zeta=
\frac{u+\mu_k-k+1}{u+\mu_k-k+2}\ts s_{m+l,k}T_{aa}(u)\ts\zeta'
\non\eea
for $a>m+l$, and
\bea\label{Taamu2}
T_{aa}(u)\zeta=
s_{m+l,k}T_{aa}(u)\ts\zeta'
\non\eea
for $a=m+l$. If $a<m+l$ then (\ref{Ttaumu2}) with $a=b$ gives
\bea\label{Taamu3}
T_{aa}(u)\zeta=
\frac{u+\mu_k-k}{u+\mu_k-k+1}\ts s_{m+l,k}T_{aa}(u)\ts\zeta'.
\non\eea
In all the cases the relation (\ref{Taazeta}) 
follows by an obvious induction.

Let us now prove (\ref{szeta}). Suppose that $\mu_k=\lambda_{k+l}$
but $\wt \zeta:=s_{m+l,k}\ts \zeta\ne 0$. 
Then $\wt \zeta$ is a $\gl(m)$-highest vector of weight 
$\wt\mu=\mu-\delta_k$. This is only possible if $\mu_k>\mu_{k+1}$
which will be assumed.
We can repeat the previous arguments where the vector $\zeta$
is replaced with $\wt \zeta$ to show that
$\wt \zeta$ is annihilated by $T_{ab}(u)$
with $a<b$ and $\wt \zeta$ is an eigenvector for the $T_{aa}(u)$.
This means that $\wt \zeta$ is the highest vector of the $\Y(n)$-module
$L(\lambda)^+_{\wt\mu}$. Let us verify that
this vector is orthogonal to all elements of $L(\lambda)^+_{\wt\mu}$.
Indeed, by the Poincar\'e--Birkhoff--Witt theorem for the algebra
$\Y(n)$ (see e.g. \cite[Corollary 1.23]{mno:yc}),
$L(\lambda)^+_{\wt\mu}=\Y_-\ts \wt \zeta$, where $\Y_-$
is the span of monomials in the coefficients of the series $T_{ba}(u)$
with $b>a$. However, due to (\ref{T*}) we have for any 
$\eta\in L(\lambda)^+_{\wt\mu}$:
\bea\label{formz}
\langle \wt\zeta,T_{ba}(u)\eta\rangle=
\langle T_{ab}(u)\wt\zeta,\eta\rangle=0.
\non\eea
It remains to check that
$\langle \wt\zeta,\wt\zeta\rangle =0$. Indeed, by (\ref{conj})
\bea\label{wtzeta}
\langle \wt\zeta,\wt\zeta\rangle={\rm const}\cdot
\langle \zeta,s_{k,m+l}s_{m+l,k}\zeta\rangle.
\non\eea
By (\ref{psiTab2}) we have
\beq\label{sT}
s_{k,m+l}s_{m+l,k}\zeta=(-1)^m T_{m+l,m+l}(-h_k+1)\ts \zeta.
\eeq
Note that 
\bea\label{hk}
h_k\ts \zeta=(\mu_k-k+1)\ts \zeta=(\lambda_{k+l}-k+1)\ts \zeta.
\non\eea
Therefore, (\ref{sT}) equals 
$(-1)^m T_{m+l,m+l}(-\lambda_{k+l}+k)\ts \zeta$ which is zero
by (\ref{Taazeta}). Thus, $\wt\zeta$ has to be zero
which contradicts to our assumption.

Finally, to show that $\zeta\ne 0$ we apply appropriate operators $s_{ia}$
to $\zeta$ repeatedly to get the highest vector $\xi$ of $L(\lambda)$
with a nonzero coefficient. Indeed, by (\ref{psiTab2}) we have
\bea\label{sz}
s_{k,m+l}\ts\zeta=(-1)^m T_{m+l,m+l}(-\mu_k+k-1)\ts\zeta'
\non\eea
which equals 
\beq\label{Tzeta}
-(\mu_k-\nu_{l}^{(1)}-k+1)\cdots (\mu_k-\nu_{l}^{(m+1)}-k+m+1)\ts\zeta'
\eeq
by (\ref{Taazeta}). Here the numbers $\nu_{l}^{(i)}$ are defined by
(\ref{nu}) with 
\bea\label{mu}
\mu=(\lambda_1,\dots,\lambda_{k-1},\mu_k+1,\mu_{k+1},\dots,\mu_m).
\non\eea
So, the coefficient at $\zeta'$ in (\ref{Tzeta})
is nonzero. \endproof

Theorem \ref{thm:hvpr} implies the following
corollary.
As usual,
by the {\it content\/} of a cell $\alpha=(i,j)\in\lambda/\mu$ we mean
the number $j-i$.

\bco 
The Drinfeld polynomials for
the $\Y(n)$-module $L(\lambda)^+_{\mu}$ are given by
\beq\label{lamdr}
P_{a}(u)=\prod_c (u+c),\qquad a=1,\dots,n-1,
\eeq
where $c$ runs over the contents of the top cells of columns
of height $a$ in the diagram $\lambda/\mu$. 
\eco

\Proof By the definition (\ref{drpol}) and the formula (\ref{lambda})
we have
\bea\label{drlam}
P_a(u)=\prod_{k=1}^{m+1} (u+\nu_{a+1}^{(k)}-k+1)\cdots 
(u+\nu_{a}^{(k)}-k),
\non\eea
where the $k$th factor is assumed to be equal to $1$ if 
$\nu_{a}^{(k)}=\nu_{a+1}^{(k)}$.
Now (\ref{nu}) implies that this product
coincides with (\ref{lamdr}). \endproof

If $\lambda=(10,8,5,4,2)$
and $\mu=(6,3)$ (see the example above) then we have
\bea\label{Dp}
P_{1}(u)=(u+4)(u+8)(u+9),\qquad P_2(u)=u(u+3)(u+6)(u+7).
\non\eea

The corollary shows that the $\Y(\sll(n))$-module $L(\lambda)^+_{\mu}$
is isomorphic to the corresponding
elementary representation of $\Y(\sll(n))$
which was studied in~\cite{nt:ry}.

\begin{\bib}{99}

\bibitem{ast:po}
{R. M. Asherova, Yu. F. Smirnov and V. N. Tolstoy},
{\it Projection operators for simple Lie groups},
Theor. Math. Phys. {\bf 8} (1971), 813--825.

\bibitem{cp:yr}{V. Chari {\rm and} A. Pressley},
{\it Yangians and $R$-matrices}, {L'Enseign. Math.} {\bf
36} (1990), 267--302.

\bibitem{c:ni}
{I. V. Cherednik},
{\it A new interpretation of Gelfand--Tzetlin bases}, {Duke Math. J.}
{\bf 54}
(1987),
563--577.

\bibitem{d:ae} J. Dixmier, {\it Alg\`ebres enveloppantes\/},
Gauthier-Villars, Paris, 1974.

\bibitem{d:ha}
{V. G. Drinfeld},
{\it Hopf algebras and the
quantum Yang--Baxter equation}, 
{Soviet Math. Dokl.} {\bf
32} (1985), 254--258.

\bibitem{d:nr}
{V. G. Drinfeld},
{\it A new realization of
Yangians and quantized affine algebras}, {Soviet Math. Dokl.} {\bf
36} (1988), 212--216.

\bibitem{gkllrt:ns} 
{I. M. Gelfand, D. Krob, A. Lascoux,
B. Leclerc, V. S. Retakh and J.-Y. Thibon},
{\it Noncommutative symmetric functions}, Adv. Math. 
{\bf 112} (1995), 218--348.

\bibitem{gr:dm} {I. M. Gelfand and V. S. Retakh}, {\it Determinants
of matrices over noncommutative rings}, {Funct. Anal. Appl.} {\bf 25}
(1991), 91--102.

\bibitem{gs:ci}
{M. D. Gould and N. I. Stoilova}, {{\it Casimir invariants and
characteristic identities for} ${\gl}(\infty)$},
{J. Math. Phys.} {\bf 38} (1997), 4783--4793.

\bibitem{kl:mi}
{D. Krob and B. Leclerc}, {\it Minor identities for 
quasi-determinants and quantum determinants}, {Comm. Math. Phys.}
{\bf 169} (1995), 1--23.

\bibitem{krs:yb}
{P. P. Kulish, N. Yu. Reshetikhin and E. K. Sklyanin},
{\it Yang--Baxter equation and representation theory}, {Lett. Math. Phys.}
{\bf 5}
(1981),
393--403.

\bibitem{ks:qs}
{P. P. Kulish and E. K. Sklyanin},
{\it Quantum spectral transform method: recent developments},
in \lq Integrable Quantum Field Theories', {Lecture Notes in Phys.}
{\bf 151}
Springer,
Berlin-Heidelberg,
1982,
pp. 61--119.

\bibitem{m:fd} {A. I. Molev},
{\it Finite-dimensional irreducible representations of twisted
Yangians},
{J. Math. Phys.} {\bf 39} (1998), 5559--5600.

\bibitem{mno:yc}
{A. Molev, M. Nazarov and G. Olshanski},
{\it Yangians and classical Lie algebras}, 
Russian Math. Surveys
{\bf 51}:2
(1996),
205--282.

\bibitem{mo:cc}
{A. Molev and G. Olshanski},
{\it Centralizer construction for twisted Yangians},
Preprint CMA 065-97, Austral. Nat. University, Canberra;
q-alg/9712050.

\bibitem{nt:yg}
{M. Nazarov and V. Tarasov},
{\it Yangians and Gelfand--Zetlin bases}, Publ. RIMS, Kyoto Univ.
{\bf 30} (1994), 459--478.

\bibitem{nt:ry}
{M. Nazarov and V. Tarasov},
{\it Representations of Yangians with Gelfand--Zetlin bases}, 
{J. Reine Angew. Math.} {\bf 496} (1998), 181--212.

\bibitem{o:qi}
{A. Okounkov},
{\it Quantum immanants and higher Capelli identities},
{Transformation Groups}
{\bf 1}
(1996),
99--126.

\bibitem{oo:ss}
A. Okounkov and G. Olshanski,
{\it Shifted Schur functions},
St.\,Petersburg Math. J.
{\bf 9}
(1998), 239--300.

\bibitem{o:ea} {G. I. Olshanski}, 
{\it Extension of the algebra $U(g)$ for 
infinite-dimensional classical Lie algebras $g$, 
and the Yangians $Y(gl(m))$.}
{Soviet Math. Dokl.} {\bf 36} (1988), no. 3, 569--573.

\bibitem{o:ri}{G. I. Olshanski},
{\it Representations of 
infinite-dimensional classical groups, limits of enveloping algebras, and
Yangians}, in \lq Topics in Representation Theory' (A.~A.~Kirillov, Ed.),
{Advances in Soviet Math.} {\bf 2}, AMS, Providence RI, 1991, pp.
1--66.

\bibitem{o:ty}
{G. I. Olshanski},
{\it Twisted Yangians and infinite-dimensional classical Lie algebras},
in \lq Quantum Groups' (P.~P.~Kulish, Ed.), {Lecture Notes in Math.}
{\bf 1510},
Springer,
Berlin-Heidelberg,
1992,
pp. 103--120.

\bibitem{tf:qi}{L. A. Takhtajan {\rm and} L. D. Faddeev},
{\it Quantum inverse scattering method and the Heisenberg $XYZ$-model},
{Russian
Math. Surv.} {\bf 34} (1979), no. 5, 11--68.

\bibitem{w:cg}
{H. Weyl},
{\it Classical Groups, their Invariants and Representations},
{Princeton Univ. Press},
Princeton NJ,
1946.

\bibitem{z:cl}
{D. P. \v{Z}elobenko},
{\it Compact Lie groups and their representations},
{Transl. of Math. Monographs}
{\bf 40} AMS,
Providence RI
1973.

\bibitem{z:sa} 
{D. P. Zhelobenko},
{\it $S$-algebras and Verma modules
   over reductive Lie algebras}, 
Soviet. Math. Dokl.
{\bf 28} (1983), 696--700.

\bibitem{z:za} 
{D. P. Zhelobenko},
{\it $Z$-algebras over reductive Lie
   algebras}, Soviet. Math. Dokl.
{\bf 28} (1983), 777--781.

\bibitem{z:it} 
{D. P. Zhelobenko},
{\it An introduction to the theory of $S$-algebras over reductive
Lie algebras}, in \lq Representations of
Lie groups and Related Topics' 
(A.~M.~Vershik, D.~P.~Zhelobenko,~Eds.), {Adv. Studies in Contemp. Math.}
{\bf 7},
Gordon and Breach Science Publishers,
New York,
1990,
pp. 155--221.

\end{\bib}

\end{document}